\documentstyle[12pt]{article}
\newcommand{\myref}[1]{\hangafter=1\hangindent=2pc\parindent= 0pt
     \baselineskip=12pt #1  \parskip=10pt\par\baselineskip=20pt}
\newcommand{\myleft}[1]{\dimen0=\parindent\hangafter=1
     \hangindent=2pc\parindent=0pt #1  \parskip=10pt\par\parindent=\dimen0}
\newcommand{\Ho}{{\rm H}_o}
\newcommand{\dconv}{{\,\stackrel{2\,}{\Rightarrow}}}
\newcommand{\mtx}{m^{\rm T}X}
\newcommand{\mtPx}{m^{\rm T}{\bf P}x}
\newcommand{\mbar}{{\overline{m}}}
\newcommand{\Mbar}{{\overline{M}}}
\newcommand{\xbar}{{\overline{x}}}
\newcommand{\Xbar}{{\overline{X}}}
\newcommand{\Fhat}{{\hat F}}
\newcommand{\mdt}{(m - \mbar )^{\rm T}}
\newcommand{\mdiff}{(m_i - \mbar )}
\newcommand{\law}{{\cal L}}
\newcommand{\mdtPx}{(m-\mbar)^{\rm T}{\bf P}x}
\newcommand{\norm}[1]{{\|#1\|}}
\newcommand{\Pbf}{{\bf P}}
\newcommand{\one}{{\bf 1}}
\newcommand{\E}{{\rm E}}
\newcommand{\Var}{{\rm Var}}
\newcommand{\M}{{\cal M}}

\newcommand{\Lbar}{{\overline{L}}}
\newcommand{\Hajek}{H{\'a}jek}
\newcommand{\Sidak}{{\v S}id{\'a}k}
\newcommand{\Cibisov}{{\v C}ibisov}
\newcommand{\reals}{{\cal R}}
\newtheorem{theorem}{Theorem}
\newtheorem{lemma}{Lemma}
\newtheorem{corollary}{Corollary}

\title{Inefficient Best Invariant Tests}
\author{Richard A. Lockhart\thanks{This manuscript was written in 1992 but the
reviewers were so caustic I never tried again to publish it.  It seems relevant today.
Supported by the 
Natural Sciences and Engineering
Research Council of Canada.}\\
Department of Mathematics and Statistics\\
Simon Fraser University\\
Burnaby, BC\\
Canada V5A 1S6}
\begin{document}
\maketitle
\begin{abstract}
Test statistics which are invariant under various subgroups of the
orthogonal group are shown to provide tests whose powers are 
asymptotically equal to their level 
against the usual type of contiguous alternative in models
where the number of parameters is allowed to grow as the 
sample size
increases. The result is applied to the usual analysis of variance test 
in the Neyman-Scott many means problem and to an analogous 
problem in exponential families.  Proofs are based on a method 
used by \Cibisov(1961)
to study spacings statistics in a goodness-of-fit problem.  We
review the scope of the technique in this context.
\end{abstract}
\myleft{{\em Keywords.}
Asymptotic relative efficiency, Neyman-Scott many means problem,
goodness-of-fit, spacings statistics, many parameter problems, 
permutation
central limit theorem, bootstrap.}
\myleft{{\em AMS 1980 Classification.} Primary 62G20, Secondary
62G30.} 
\myleft{Running title: Invariant tests}

\section{Introduction}

\par

Consider the problem of data $X$ coming from a model indexed by a 
parameter
space $\cal M$.  Suppose there is a group $G$ which acts both on 
the data and on the parameter space so that $gX$ has the 
same distribution under $m\in \cal M$
as $X$ has under $g^{-1}m.$  
The problem of testing $\Ho : m\in\M_0\subset\M$ is
invariant under $G$ if $m\in\M_0$ iff  $gm\in\M_0$ for all  $g\in 
G.$ 
In 
what 
follows we shall impose the stronger condition that $gm=m$ for all 
$m\in\M_0$ and all $g\in G$.

In this note we use the observation that for an invariant test 
the power at alternative $m$ minus the level of the test is the 
covariance
between the test function and the likelihood ratio averaged over the 
orbit of
$m$ under $G$ to study the power of invariant tests under 
contiguous
alternatives.  The technique was used in \Cibisov (1961) to study the
asymptotic behaviour of tests of uniformity based on sample 
spacings.  

Our results may be summarized as follows.  When for a given point 
$m$ in the 
null hypothesis, the number of possible directions of departure from 
$m$ into the alternative hypothesis 
grows with the sample size the power of invariant tests may be
expected to be low.  We will exhibit a variety of examples in which
the power minus the level of invariant tests converges to 0 
uniformly in the class of invariant tests.

We begin in section 2 with a simple normal example to illustrate 
the technique and set down the basic identities and inequalities.  
In section 3 we extend the normal example to
general exponential families using a version of the permutation 
central limit theorem.  In section 4 we examine the Neyman Scott 
many means 
problem and
extend our results to more general models.  In section 5 we revisit 
\Cibisov's example to show that a variety of
goodness-of-fit tests have ARE 0 over a large class of alternatives.  
The results of the first 5 sections suggest that either invariance is not 
always desirable or that contiguity calculations are not always the 
right way 
to compare powers of tests in such models.
Section 6 is a discussion of the relevance of contiguity calculations in 
this 
context, together with some open problems and suggestions for 
further work. An appendix contains some technical steps in the proofs.

\section{Basic Results: A Normal Example}

Suppose $X\in\reals^n$ has a multivariate normal distribution with
mean vector $m\in\reals^n$ and identity covariance matrix.  The
problem of testing the null hypothesis $\Ho\colon m=0$ 
is invariant under
the group of orthogonal transformations.  Under $\Ho$, $X$ and 
$\Pbf X$
have the same distribution for any orthogonal matrix $\Pbf$.  
Consider a
simple alternative $m \ne 0$.  The likelihood ratio of $m$ to 0 is 
$L(X) = 
\exp\{\mtx-\norm{m}^2/2\}$. Thus the Neyman Pearson test rejects 
when
$\mtx / \norm{m}$ is too large.  Under $\Ho$ this statistic has a
standard normal distribution while under the alternative $m$ the 
mean
is shifted to $\norm{m}$. Thus as $n\to\infty$ non-trivial  limiting 
power (i.e. a limiting power larger than the level and less than 1) 
results when $\norm{m}\to\delta\ne 0$.  (Throughout this paper
objects named by Roman letters depend on $n$; wherever possible 
the
dependence is suppressed in the notation.  Objects named by Greek
letters do not depend on $n$.)

We may analyse the efficiency of an invariant test relative to the
Neyman Pearson test for a given sequence of alternatives as follows.
Let ${\cal T}$ be the class of all test functions $T(X)$ for testing
$\Ho$ which are invariant under the group of orthogonal
transformations, that is, for which $T(\Pbf X) = T(X)$ for each
orthogonal transformation $\Pbf$.  We have the following theorem.

\begin{theorem}{ As $n\to\infty$,
\[\sup\{|\E_m(T)-\E_0(T)|: 
T\in{\cal T}, 
\norm{m}\le\delta\}\to 0.\]}
\end{theorem}

Since the parameter space depends on the sample size the usual 
definition of 
relative efficiency does not make sense.  Instead, given an 
alternative $m$ 
we will define the efficiency of a (level $\alpha$ where $\alpha$ is
fixed) test relative to the Neyman Pearson test 
to be $1/c^2$ where $c$ is chosen so that the test under 
consideration has 
power against $cm$ equal to the Neyman Pearson power against the 
alternative $m$.  
In traditional finite dimensional parametric models this notion 
agrees with 
the usual notion of Pitman relative efficiency (asymptotically).  Thus 
ARE will be the limit of $1/c^2$; the alternative sequence $m$ must 
be 
contiguous so that the asymptotic power of the Neyman Pearson test 
is not 1. For a discussion of the relevance of the Neyman-Pearson test as
a standard for the relative efficiency see section 6.

\begin{corollary}
{The best invariant test of $\Ho$ has ARE 0.}
\end{corollary}

\leftline{\em Proof of Theorem 1}

Let $L(X) $ be the likelihood ratio for $m$ to $0$. 
Then\[\E_m\left(T(X)\right)=\E_0\left(T(X)L(X)\right).\]
Since $\Pbf X$ and $X$ have the same distribution under $\Ho$ we 
have
\begin{eqnarray}
\E_m\left(T(X)\right)&=&\E_0\left(T(\Pbf X)L(\Pbf 
X)\right)\nonumber\\
&=&\E_0\left(T(X)L(\Pbf X)\right)
\end{eqnarray}
for all orthogonal $\Pbf.$  
Since $\Pbf$ appears only on the right
hand side of $(1)$ we may average over orthogonal $\Pbf$ to obtain
\begin{eqnarray}
\E_m\left(T(X)\right)&=&\int\E_0\left(T(X)L(\Pbf X)\right) F(d\Pbf 
)\nonumber\\
&=&\E_0\left(T(X)\int L(\Pbf X) F(d\Pbf )\right)
\end{eqnarray}
where $F$ is any probability measure on the 
compact group  of orthogonal matrices. Let 
\[\Lbar (X)=\int L(\Pbf X) F(d\Pbf ).\]  Then
\begin{eqnarray}
|\E_m(T)-\E_0(T)| &=&|\E_0(T(\Lbar -1))|\nonumber\\
&\le&\E_0(|\Lbar -1|).
\end{eqnarray}
Since the last quantity 
is free of $T$ we need only show that 
\[\sup\{\E_0(|\Lbar -1|) : \norm{m}\le\delta\}\to 0.\]
Since $\E(\Lbar)=1$ the dominated convergence theorem shows
that it suffices to prove for an arbitrary 
sequence of alternatives $m$ with $\norm{m}\le\delta$
and for a suitably chosen sequence of measures $F$ 
that $\Lbar \to 1$ in probability.  We take $F$ to be Haar
measure on the compact group of orthogonal transformations; that is, 
we give $\Pbf $ a uniform distribution.

For each fixed $X$, when $\Pbf$ has the distribution for which
$F$ is Haar measure on the orthogonal group, the vector $\Pbf X$ has 
the uniform distribution on the sphere
of radius $\norm{X}.$  Using the fact that a standard multivariate 
normal
vector divided by its length is also uniform on a sphere we find that 
\[\Lbar = H(\norm{m}\cdot\norm{X})/H(0)\]
where\[H(t) = \int_0^\pi\exp(t\cos\theta) \sin^{n-2}\theta 
\,d\theta.\]
Standard asymptotic expansions of Bessel functions 
(see, e.g. Abramowitz and Stegun, 1965, p 376ff) then make it
easy to show that $\Lbar\to 1$ in probability, finishing the
proof.$\bullet$

A test invariant under the group of orthogonal 
transformations is a function of $\norm{X}^2$ and the 
analysis above can be made directly and easily 
using the fact that this statistic has a chi-squared 
distribution with non-centrality parameter $\norm{m}^2$.  Our
interest centres on the technique of proof.  Equations (1-3) and
the argument following (3) use only the group structure of
the problem and the absolute 
continuity of the alternative with respect the null.  The remainder of 
the argument depends on an asymptotic approximation to the 
likelihood
ratio averaged over alternatives.  Whenever such an approximation 
is available we can expect to obtain efficiency results for the family 
of
invariant tests.  In the next section we apply the technique to a more 
general model by replacing the explicit Bessel function calculation 
with a 
version of the permutation central limit theorem.

\section{Exponential Families}

Suppose now that $X = (X_1,\ldots,X_n)^{\rm T}$ with the $X_i$ 
independent and $X_i$ having the exponential family density $\exp( 
m_i
x_i -\beta(m_i))$ relative to some fixed measure.  We assume the
$m_i$'s take values in $\Theta$ an open subset of $\reals.$ The
parameter space is then $\M = \Theta^n.$ Let $\Theta_0$ be a fixed 
compact subset of $\Theta.$

Consider the null hypothesis $\Ho\colon  m_1=\cdots=m_n.$ This 
problem is
invariant under the subgroup of the orthogonal group consisting of 
all permutation matrices $\Pbf .$ Let $\mbar = \sum m_i /n$; we 
also 
use $\mbar$ to denote the vector of length $n$ all of whose entries 
are
$\mbar $. The calculations leading to (1-3) establish that for any
test $T$ in $\cal T$, the family of all tests invariant under
permutations of the entries of $X,$ we have
\[|\E_m(T)-
\E_{\mbar}(T)|
\le \E_{\mbar}(|\Lbar -1|).\] where  now\[\Lbar
= \exp\left\{-\sum \left(\beta(m_i)- \beta (\mbar ) 
\right)\right\}\sum_\Pbf
\exp\left(\mdt\Pbf X\right)/n!\] 
is the likelihood ratio averaged over all
permutations of the alternative vector $m.$

\subsection{Heuristic Computations}

Think now of $X$ as fixed and $\Pbf$ as a randomly chosen 
permutation matrix. Then $\mdt\Pbf X = \mdt\Pbf (X - {\Xbar})$ has 
moment 
generating function
\[g(s)= \sum_{\Pbf}\exp\{s\mdt\Pbf (X - {\Xbar})\}/n!.\] 
If $\max\{|m_i-\mbar|: 1\le
i \le n\}\to 0$ the permutation central limit theorem suggests that
$\mdt\Pbf X$ has approximately a normal distribution with mean
0 and variance $\sum\mdiff^2\sum(X_i-{\Xbar})^2/n$.  Thus
heuristically\[g(s)\approx \exp \{ s^2 \sum \mdiff^2
\sum(X_i-{\Xbar})^2/(2n)\}.\] 
Under the same condition, $\max\{|m_i-\mbar|: 1\le i \le
n\}\to 0$, we may expand\[\sum ( \beta ( m_i) - \beta ({\mbar})) 
\approx 
\sum
\mdiff^2\beta''({\mbar})/2.\] We are led to the heuristic calculation
\[\Lbar(X)\approx\exp\{\sum\mdiff^2(S^2-\beta''({\mbar}))/2\}\] 
where
$S^2=\sum(X_i-{\Xbar})^2/n$ is the sample variance.  Since ${\rm
Var}_{\mbar}(X)=\beta''({\mbar})$ we see that $\Lbar$ should 
converge to 1
provided $\sum\mdiff^2 = o ( n^{1/2})$ or $\norm{m-\mbar} = 
o(n^{1/4}).$
In the simple normal example
of the previous section we are actually able to prove (again using
asymptotic expansions of Bessel functions) this
strengthened version of Theorem 1, namely, if $a=o(n^{1/4})$ then
\[\sup\{|\E_m(T)-\E_0(T)|: T\in{\cal T},\norm{m}\le  a\}\to 0.\] 

In the more general exponential family problem we do not know a 
permutation central limit theorem which extends  in a useful way
to give convergence of moment 
generating functions.  For contiguous alternative sequences we are 
able to
replace the moment generating function by a characteristic
function calculation.  We can then prove that
invariant tests have power converging to their level uniformly on
compact subsets of $\Theta.$

\begin{theorem}
As $n\to\infty$,\[\sup\{|\E_m(T)-\E_{\mbar}(T)|: 
T\in{\cal T}, \norm{m - \mbar }\le\delta,m_i\in\Theta_0\}\to 0.\]
\end{theorem}

\subsection{Proof of Theorem 2}

Our proof uses contiguity techniques to replace the moment 
generating function used in the heuristics above by a corresponding
characteristic function calculation.  A standard compactness 
argument reduces the problem to showing that
\[\E_m(T)-\E_\mbar(T)\to 0\] for an arbitrary sequence of 
alternatives
$m$ satisfying $\sum\mdiff^2\le\delta^2$, $m_i\in\Theta_0$, and an
arbitrary sequence $T\in {\cal T}$.  Our proof is now in three stages.
First we prove that any such sequence is contiguous to the null
sequence indexed by $\mbar$.  The second step is to eliminate the
particular statistics, $T$, as at (3), reducing the problem to showing
that the permutation characteristic function of the log-likelihood
ratio is nearly non-random.  The final step is to establish the latter
fact appealing to Theorem 3.

\vskip\parskip
\leftline{{\em Step} 1: Contiguity of $m$ to $\mbar$.}

The log-likelihood ratio of $m$ to
$\mbar$ is\[\ell(X)= \sum \mdiff (X_i-\beta({\mbar}))-\sum
(\beta(m_i)-\beta({\mbar})).\]
Under $\mbar$ the mean of $\ell(X)$ is 
\[-\sum
(\beta(m_i)-\beta({\mbar}))=-\sum(m_i-\mbar)^2\beta''(t_i)\]
for some $t_i$ between $m_i$ and $\mbar$.  In view of the 
conditions on $m$ and the compactness of $\Theta_0$ this
sequence of means is bounded.  Also under $\mbar$ the variance of 
$\ell(x)$ is \[\sum (m_i-\mbar)^2\beta''(\mbar).\]  Again this is
bounded.  Thus the sequence of log-likelihood ratios is tight under
the null sequence $\mbar$ and so the sequence of alternatives, $m$, 
is
contiguous to the null sequence $\mbar$.
 
\begin{lemma}
 Suppose $Q$ is a sequence of measures contiguous to a sequence of 
 measures $P$.
 If $T$ is a bounded sequence of statistics such that 
 \[{\rm E}_P(T\exp(i\tau\log dQ/dP))-
 {\rm E}_P(T){\rm E}_P(\exp(i\tau\log dQ/dP))\to 0\]
 for each real $\tau$ then \[{\rm E}_Q(T)-{\rm E}_P(T)\to 0.\]
\end{lemma}

\myleft{{\em Remark}: In the Lemma the random variable $\log 
dQ/dP$ 
 can be replaced by any random variable $S$ such that $\log dQ/dP-
S$ tends 
 to 0 in probability under $P$.}

\myleft{{\em Remark}: The Lemma is very closely connected with 
LeCam's
third lemma (see \Hajek\null\  and \Sidak, 1967, page 209, their 
formula 4) which
could also be applied here to the $Q$ characteristic function of $T$.}

Before proving the lemma we finish the theorem.

\vskip\parskip
\leftline{{\em Step} 2: Elimination of the sequence $T$.}

Let $H(t,X)=\exp(it\ell(X))$. According to the lemma we must show
\begin{equation}
 \label{four}
 \E_\mbar(T(X)H(\tau,X))-\E_\mbar(T)\E_\mbar(H(\tau,X)) \to 0\, .
\end{equation}
Arguing as in (1-3) the quantity in (\ref{four}) may be seen to be
\[\E_\mbar(T(X){\overline{ H}}(\tau,X)) - 
 \E_\mbar(T(X))\E_\mbar({\overline{H}}(\tau,X))\] 
where ${\overline{H}}(\tau,X)=\sum_\Pbf H(\tau,\Pbf X)/n! = 
\E(H(\tau, \Pbf X)|X)$. 
Since both $T$ and ${\overline{H}}(\tau,X)$ 
are bounded it suffices to prove that 
\begin{equation}
 \label{hbarcond}
 {\overline{H}}(\tau,X)-\E_\mbar({\overline{H}}(\tau,X)) \to 0
\end{equation} 
in probability.

\myleft{{\em Remark}: If $\ell(X)=S(X)+o_P(1)$ under $\mbar$ then 
in view of the remark following the lemma the 
random variable $\ell(X)$ can be replaced  in the definition of $H$ 
by the random variable $S(X)$.
Moreover, if $S^*(X)=S(X)+a$ then (\ref{hbarcond}) 
will hold for $S^*$ replacing 
$\ell$ if and only if it holds for $S$ replacing $\ell$ 
because the sequence $\exp(i\tau a)$ is bounded in modulus.}

\vskip\parskip
\leftline{{\em Step} 3: Application of the Permutation Central Limit
 Theorem}

Using the last remark take $S^*(X)=\mdtPx$.  The variable 
$\overline{H}$ 
then becomes $\sum_\Pbf\exp(i\tau\mdtPx)/n!$ whereas
$\E_\mbar({\overline{H}}(\tau,X))$ becomes
$E_\mbar(\exp(i\tau\mdt X))$.   

Now let $\hat F$ be the empirical 
distribution of the $X_1,\ldots,X_n$.
Let $X_1^*,\ldots,X_n^*$ be independent and identically distributed
according to $\hat F$. We will show in the next section that 
\begin{equation}
 \label{woutrep}
 \sum_\Pbf\exp(i\tau\mtPx)/n!-\E_{\hat F}
 \left(\exp(i\tau\sum m_jX_j^*)\right)\to 0
\end{equation}
and  
\begin{equation}
 \label{empiricalcf}
 \E_{\hat F}
 \left(\exp(i\tau\sum m_jX_j^*)\right)-
 \E_\mbar\left(\exp(i\tau\sum m_jX_j)\right)\to 0
\end{equation}
in probability for each fixed $\tau$. This will finish the proof of
Theorem 2 except for establishing the lemma.

To prove the (undoubtedly well-known) lemma
argue as follows.  Letting $\ell=\log dQ/dP$, the condition shows that
$\E_P(T\phi(\ell))-\E_P(T)\E_P(\phi(\ell))\to 0$ for each bounded 
continuous function $\phi$.
There is then a sequence $a$ tending to infinity so slowly that 
$\E_P(Tf(\ell))-\E_P(T)\E_P(f(\ell))\to 0$ where $f(x)=\min(e^x,a)$.
In view of contiguity $\E_Q(T)-\E_P(Te^\ell)\to 0$ and for any 
sequence $a$ tending
to infinity $|E_P(T(e^\ell-f(\ell)))|\le \E_P(|e^\ell-f(\ell)|)\to 0$.  The
lemma follows.$\bullet$

\subsection{The Permutation Limit Theorem}

Suppose $m \in \reals^n$ and $x\in\reals^n$ are two (non-random) 
vectors with $\mbar=0$ for convenience. Suppose $\Pbf$
is a random $n\times n$ permutation matrix. The random variable
$m^{\rm T}\Pbf x$ has mean $n\mbar\,\xbar$ and variance
$s^2=\sum\mdiff^2\sum(x_i-\xbar)^2/(n-1)$.  
An alternative description of the random variable $m^{\rm T}\Pbf x$ 
is as follows.  Let
$J_1,\ldots,J_n$ be a simple random sample drawn without 
replacement from
the set $\{1,\ldots,n\}$.  Then $J_1,\ldots,J_n$ is a random 
permutation
of $\{1,\ldots,n\}$ and $m^{\rm T}\Pbf x$ has the same distribution 
as
$\sum m_ix({J_i})$; we now use functional rather than subscript 
notation for legibility.  \Hajek (1961, see his formula 3.11) shows 
that it is possible to
construct, on a single probability space, $J_1,\ldots,J_n$ together 
with 
$J_1^*,\ldots,J_n^*$, a random sample 
drawn {\em with} replacement from $\{1,\ldots,n\}$  in such a way 
that
\begin{equation}
 \E[(\sum m_i x(J_i) - \sum m_i x(J_i^*))^2] \le 3s\,\max | x_i-\xbar| 
/
 (n-1)^{1/2}.
\end{equation}
Since the $X_i$ have an exponential family distribution and $\Theta$ 
is
compact it is straightforward to check that 
$$S\,\max|X_i-\Xbar|  /(n-1)^{1/2} =O_P(\log n/n^{1/2})=o_P(1)$$
under $\mbar$ where $S$ denotes the sample standard deviation. 
In view of the elementary inequality
\begin{equation}
 | \E(exp(itW))-\E(\exp(itW')) | \le t^2 \E((W-W')^2) +
 |t|\E^{1/2}((W-W')^2)
\end{equation}
this establishes (\ref{woutrep}).

It seems easiest to deal with (\ref{empiricalcf}) after 
recalling some facts about weak
convergence in the presence of moment conditions.  Consider the set,
$\Delta$, of distribution functions $F$ on $\reals^\pi$ which have 
finite variance.  Throughout this paper we write $F\dconv\Phi$ if
$\E_F\left(\gamma(X)\right)\to\E_\Phi\left(\gamma(X)\right)$ 
for each fixed continuous function
$\gamma$ such that $\gamma(x)/(1+\norm{x}^2)$ is bounded.  (The 
notation is
$\E_F\left(\gamma(X)\right)=\int\gamma(x)F(dx)$.) This notion of 
convergence
defines a topology on $\Delta $ which can be metrized by a metric 
$\rho_2$ in
such a way that $\Delta$ becomes a complete separable metric space. 
In fact the metric
$\rho_2$ may be taken to be the Wasserstein metric;  
see Shorack and Wellner (1986, pp 62-65). 

Suppose that $\Delta_0$ is a subset of $\Delta $. The following are 
equivalent:
\renewcommand{\theenumi}{\roman{enumi}}
\begin{enumerate}
\item $\Delta_0$  has compact closure.
\item
for each $\epsilon > 0$
there is a fixed compact subset $\Psi$ of $\reals^\pi$ such that
$\E_F(\norm{X}^21(X\not\in\Psi))\le\epsilon$ for all $F\in\Delta_0$.  
\item
there is a fixed function $\Psi$ such that
$\E_F\left(\norm{X}^21(\norm{X}\ge t)\right)\le\Psi(t)$ for all 
$F\in\Delta_0$ 
and all $t$ where $\Psi$ has the property that
$\lim_{t\to\infty}\Psi(t) = 0$.
\item
the family $\Delta_0$ makes $\norm{X}^2$ uniformly integrable.
\end{enumerate}
Notice that for
each fixed $\Psi$ such that $\lim_{t\to\infty}\Psi(t) = 0$ the family 
of distributions $F$ with
$\E_F\left(\norm{X}^21(|X|\ge t)\right)\le\Psi(t)$ for all $t$ is a 
compact 
metric space.  
Finally note that $F\dconv\Phi$ if and only if  $F$ converges in
distribution to $\Phi$ and 
$\E_F(\norm{X}^2)\to\E_\Phi(\norm{X}^2)$.  
All the foregoing results are trivial modifications of the usual results
for convergence in distribution; see Billingsley (1968) pp 31-41. It 
will be useful to let $\rho_0$ denote a metric analogous to $\rho_2$ 
for 
which the space of all distributions on $\reals^\pi$ becomes a 
complete 
separable metric space with the topology of convergence in 
distribution.  If
$\Delta_0$ is a compact subset of $\Delta$ for the metric $\rho_2$ 
then it is a compact subset of the space of all distributions for the 
metric
$\rho_0$.

To state the result let $\hat F$ be the empirical 
distribution of the numbers $X_1,\ldots,X_n$.
Let $X_1^*,\ldots,X_n^*$ be independent and identically distributed
according to $\hat F$. 
Let $\law_{\hat F}$ be the (conditional given $X_1,\ldots,X_n$) 
law of $\sum m_i X_i^*$ and $\law_F$ be the law of $\sum m_i X_i$.

%\begin{corollary}
%Under the conditions of Theorem 3, for each $\tau$,
%$$\sum_\Pbf\exp(i\tau\mtPx)/n!-\E_{\hat F}
%\left(\exp(i\tau\sum m_jX_j)\right)\to 0.$$
%\end{corollary}

%We will be using the theorem on a vector $x$ which is in fact  a 
%sample 
%from some distribution $F$. 

\begin{theorem}
\label{thm4}
Let $\Delta_0$ be a fixed compact subset of $\Delta$. Suppose $F$ is 
any sequence in $\Delta_0$ and that for each $n$ the vector
$X=(X_1,\ldots,X_n)^{\rm T}$ has independent entries distributed 
according to $F$.  Assume $m$ is an arbitrary sequence satisfying 
$\mbar=0$ and $\sum m_i^2 \le\delta^2$.
If $\hat F$ is the empirical distribution function of $X_1,\ldots,X_n$ 
then $\rho_2(\law_F,\law_{\hat F})\to 0$ in probability.
\end{theorem}

\begin{corollary}
Under the conditions of Theorem \ref{thm4}
$$\E_{\hat F}
\left(\exp(i\tau\sum m_jX_j)\right)-
\E_F\left(\exp(i\tau\sum m_jX_j)\right)\to 0$$
in probability for each fixed $\tau$.
\end{corollary}

It is straightforward to rephrase a (weakened) version of \Hajek's 
permutation central limit theorem in the notation of Theorem 
\ref{thm4}. 
Let $\law_\Pbf$ be the law of $ m^{\rm T}\Pbf x$. 

\begin{theorem}
  \label{thm3}
  Suppose that $m\in\reals^n$ and $x\in\reals^n$ are two
  sequences of vectors such that 
\begin{equation}
\label{meancond}
n\mbar\,\xbar\to 0 ,
\end{equation}
\begin{equation}
    \label{mconds}
      \sum (m_i-\mbar)^2\le \delta^2
  \end{equation}
  for a fixed $\delta$  and  such that $x$ satisfies 
  \begin{equation}
     \label{Lindeberg}
  \sum (x_i-\xbar)^21(|x_i-\xbar| \ge  t)/n\le\Psi(t)
 \end{equation}
 where $\Psi$ is a fixed function such that 
 $\lim_{t\to\infty}\Psi(t) = 0$. If $\hat F$ is the empirical 
distribution
 of the numbers $\{1,\ldots,x_n\}$ then $\rho_2(\law_\Pbf,\law_{\hat 
F})\to 0$.
\end{theorem}

The proof of Theorem \ref{thm4} is in the Appendix. 

\myleft{{\em Remark}: Theorem \ref{thm4} asserts the validity of 
the bootstrap
approximation to $L_F$.  Theorem \ref{thm3} says that the
bootstrapping carried out in Theorem \ref{thm4} by sampling from 
the list
$X_1,\ldots,X_n$ {\em with} replacement can also be carried out 
{\em without} replacement.  Notice that the result applies only to 
{\em
contrasts} in the $X_i$; the condition $\mbar=0$ is crucial to this
bootstrap interpretation of Theorem \ref{thm3}.}

\section{Extensions}

\subsection{Non-exponential families}

The exponential family model gives the log-likelihood ratio  a rather
special form.  For general models, however, a Taylor expansion 
can be used to show that the log-likelihood ratio has almost that 
form.
Rather than try to discover the weakest possible conditions on a 
general
model permitting the conclusion we try to illustrate the idea with 
assumptions which are far from best possible.

Suppose $X=(X_1,\ldots,X_n)^{\rm T}$ with the individual $X_i$ 
independent
and $X_i$ having density  $\exp(\phi(\cdot;m_i))$ where 
$m_i\in\Theta$ 
an open subset of $\reals$.  Again consider the null hypothesis
${\Ho}\colon m_1 = \cdots = m_m$. Let $\Theta_0$ denote some 
fixed 
compact subset of $\Theta$.
The log-likelihood ratio for $m$ to $\mbar$ is
$\ell(X)=\sum\{\phi(X_i;m_i)-\phi(X_i;\mbar)\}$. 
Assume that $\phi$ is twice differentiable with respect to the 
parameter.
Let $\phi_i$ denote the $i$-th derivative of $\phi$ with respect to 
the
parameter.  Assume that the usual identities of large sample theory 
of
likelihood hold, namely, that $\phi_1(X,\mbar)$ has mean 0 and 
finite
variance $\iota(\mbar)$ under $\mbar$ and that $\phi_2(X,\mbar)$ 
has 
mean $-\iota(\mbar)$.  Then we may write 
\begin{eqnarray}
\label{taylor}
\ell(X)&=&\sum(m_i-\mbar)\phi_1(X_i,\mbar)-\sum(m_i-\mbar)^2 
\iota(\mbar)/2\cr
&&\mbox{}+\sum(m_i-
\mbar)^2\{\phi_2(X_i,\mbar)+\iota(\mbar)\}/2 \cr
&&\mbox{}+ \sum(m_i-\mbar)^2 \{\phi_2(X_i,t_i)-
\phi_2(X_i,\mbar)\}
\end{eqnarray}
where $t_i$ is between $m_i$ and $\mbar$.

We see that, under the usual sort of 
regularity conditions, there will exist a constant $\alpha$ such that
\begin{equation}
\label{meanbound}
\left|\E_m(\phi(X_i;m^*)-\phi(X_i;m))\right| \le \alpha(m^*-m)^2
\end{equation}
and such that 
\begin{equation}
\label{varbound}
{\rm Var}_m(\phi(X_i;m^*)-\phi(X_i;m)) \le \alpha(m^*-m)^2
\end{equation}
for all $m$ and $m^*$ in $\Theta_0$. Under these two conditions we 
see that
$|\E(\ell(X))|\le\alpha\sum(m_i-\mbar)^2$ and ${\rm
Var}_\mbar(\ell(X))\le\alpha\sum(m_i-\mbar)^2$.  Thus any 
sequence $m$ with
$\sum(m_i-\mbar)^2\le\delta^2$ and all $m_i\in\Theta_0$ is 
contiguous 
to the null sequence $\mbar$.
The first two steps in the proof of Theorem 3 may therefore be seen 
to apply to general one-parameter models under regularity 
conditions,
specifically whenever (\ref{meanbound}) and (\ref{varbound}) hold. 

Further regularity assumptions are necessary in order for Step 3 of
Theorem 2 to go through in the present context.  To get some insight
consider the situation where $\max\{|m_i-\mbar|; 1\le i\le n\}\to 0$.
Under further regularity conditions we will have
\begin{equation}
\label{equicontinuity}
\sum(m_i-\mbar)^2\{\phi_2(X_i;t_i)-\phi_2(X_i;\mbar)\}\to 0
\end{equation}
and
\begin{equation}
\label{fisher}
\sum(m_i-\mbar)^2\{\phi_2(X_i;\mbar)+\iota(\mbar)\}\to 0
\end{equation}
in probability. Assuming that (\ref{equicontinuity}) and 
(\ref{fisher})
hold we see
\begin{equation}
\label{asylike}
\ell(X)=\sum(m_i-\mbar)\phi_1(X_i;\mbar)-\sum(m_i-\mbar)^2 
\iota(\mbar)/2
+o_P(1).
\end{equation}
Define $S^*(X) = \sum(m_i-\mbar)\phi_1(X_i;\mbar)$.  Step 3 of the 
proof
of Theorem 2 may now be carried through with $X_i$ replaced by
$U_i=\phi_1(X_i;\mbar)$ provided that the map from $\Theta_0$ to 
$\Delta$
which associates $\mbar$ with the $\mbar$ distribution of
$\phi_1(X_i;\mbar)$ is continuous.

The assumption that  
$\max\{|m_i-\mbar|; 1\le i\le n\}\to 0$ can be avoided; we now 
make this assertion precise.  We will need three more assumptions:
\begin{equation}
\label{fisherunifint}
\sup\{\E_m(|\phi_2(X,m)+\iota(m)|1(|\phi_2(X,m)+\iota(m)| \ge t)); 
m\in\Theta_0\} \to 0
\end{equation}
as $t\to\infty$. Define 
$$W(X,\epsilon,m)= \sup\{|\phi_2(X,m^\prime) - \phi_2(X,m)|; 
|m^\prime -m| \le \epsilon, m^\prime\in\Theta_0\}.$$
The second assumption will then be
\begin{equation}
\label{fishercont}
\lim_{\epsilon\to 0}\sup\{\E_m(W(X,\epsilon,m)); m\in\Theta_0\} = 
0.
\end{equation}
Finally we will assume that the map
\begin{equation}
\label{scoreassump}
(m,\mbar)\mapsto\cases{\law((\phi(X,m)-\phi(X,\mbar))/(m-
\mbar)|\mbar)&$m\ne\mbar$\cr
\law(\phi_1(X,\mbar))&$m = \mbar$}
\end{equation}
is continuous from $\Theta_0\times\Theta_0$ to $\Delta$ where 
$\law(X|m)$ denotes the 
law of $X$ when $m$ is true.
\begin{theorem}
\label{nonexpthm}
Assume conditions (\ref{meanbound}, \ref{varbound}, 
\ref{fisherunifint}, 
\ref{fishercont}, \ref{scoreassump}). Then as  $n\to\infty$, 
\[\sup\{|\E_m(T)-
\E_{\mbar}(T)|: 
T\in{\cal T}, \norm{m - \mbar }\le\delta,m_i\in\Theta_0\}\to 0.\]
\end{theorem}

\leftline{\em Proof}
Assume, without loss of generality that the entries in $m$ have been 
sorted 
so that $|m_1-\mbar| \ge \cdots \ge |m_n-\mbar|$.  Let $k=k(n)$ be 
any 
sequence tending to infinity.  Then $\sum(m_i-\mbar)^2 \le 
\delta^2$ implies 
that $\max\{|m_i-\mbar|; k\le i\le n\}\to 0$.  The assumptions now 
imply that 
\begin{eqnarray*}
\ell(X) & = & \sum_{i\le k}\{\phi(X_i,m_i)-\phi(X_i,\mbar)\} + 
\sum_{i > k}(m_i-\mbar)\phi_1(X_i,m_i) \cr & & - \sum_{i > k}(m_i-
\mbar)^2\iota(\mbar)/2 
+ o_P(1).
\end{eqnarray*}
Define 
\begin{eqnarray*}
S^*(X) & = & \sum_{i\le k}\{\phi(X_i,m_i)-\phi(X_i,\mbar) - 
\E_\mbar (\phi(X_i,m_i)-\phi(X_i,\mbar))\} \cr & & + 
\sum_{i > k}(m_i-\mbar)\phi_1(X_i,m_i).
\end{eqnarray*}

Let ${\overline{H}}(\tau,X)=\sum_\Pbf \exp(i\tau S^*(\Pbf X))/n!$. 
As 
before it suffices to prove that 
\begin{equation}
\label{step3}
{\overline{H}}(\tau,X)-\E_\mbar({\overline{H}}(\tau,X)) \to 0
\end{equation} 
in probability.  The proof may be found in the Appendix; its length is due
to our failure to impose the condition that $\max\{|m_i-\mbar|; k\le i\le
n\}\to 0$ .

\subsection{The Neyman-Scott Problem}

Consider now the Neyman-Scott many means problem in the 
following form.
Let $\{X_{ij};1\le j\le \nu ,1\le i \le n\}$ be independent normals
with mean $m_i$ and standard deviation $\sigma$.  The usual 
Analysis of Variance $F$-test  of the hypothesis that $m_1=\cdots = 
m_n$ is 
invariant under permutations of the indices $i$. 
Any level $\alpha$ test of this null hypothesis for this
model with $\sigma$ unknown is a level $\alpha$ test in any 
submodel 
with a known value of $\sigma$.  When $\sigma$ is known the 
argument of 
section 3 can be applied to the vector 
$X=({\overline{X}}_{1\cdot},\ldots,
{\overline{X}}_{n\cdot})$ of cell means to conclude that the ARE of 
ANOVA 
is 0 along any contiguous sequence of alternatives. 

This Analysis of Variance problem may be extended to the 
following multiparameter exponential family setting.  
Suppose that for $i = 1 , \ldots 
, n$  the $\reals^\pi$-valued random variable $X_{i} $ has density
$\exp\{m_i^{\rm T}x_i -\beta(m_i)\}$ relative to some fixed measure 
on
$\reals^\pi$. The natural parameter space for a single observation 
$X_i$ is some $\Theta\subset\reals^\pi$. Let $a(m)$ be some 
parameter of 
interest and consider the problem of testing $\Ho\colon a(m_1) = 
\cdots = 
a(m_n)$.  The problem is again permutation invariant.

Let $M$ be the $n\times\pi$ matrix with $i$th row $m_i$ and 
$\Mbar$ be the
$n\times\pi$ matrix with $i$th row $\mbar$. Let $X$ be the 
$n\times\pi$
matrix with $i$th row $X_i^{\rm T}$.  The log likelihood ratio of $M$ 
to
$\Mbar$ is 
$$
\ell(X)=tr((M-\Mbar)^{\rm T}X)-\sum\left(\beta(m_i)-
\beta(\mbar)\right)
$$
Denote $\norm{M}^2=tr(M^{\rm T}M)$.  Let $\Theta_0$ be some 
fixed 
compact subset of $\Theta$.  Let $\cal T$ be the family of all 
permutation
invariant test functions, $T(X)$.

\begin{theorem}
As $n\to\infty$,\[\sup\{|\E_M(T)-\E_{\Mbar}(T)|: 
T\in{\cal T}, \norm{M - \Mbar }\le\delta,m_i\in\Theta_0\}\to 0.\]
\end{theorem}

The proof of this theorem is entirely analogous to that of Theorem 2
needing only a multivariate extension of Theorems \ref{thm4} and 
\ref{thm3}.  
Suppose $M$ and $x$ are sequences of $n\times\pi$ matrices.  Let
$M_1,\ldots,M_\pi$ and $x_1,\ldots,x_\pi$ be the columns of $M$ 
and $x$ respectively.  Let $\Pbf$ be a random $n\times n$ 
permutation
matrix and let $\Fhat$ be the empirical distribution (measure on
$\reals^\pi$) of the $n$ rows of $x$.  Let $X$ be an $n\times\pi$
matrix whose rows are iid according to $\Fhat$.  Let $\law_\Pbf$ denote
the joint law of $(M_1^{\rm T}\Pbf x_1,\ldots,M_\pi^{\rm T}\Pbf x_\pi)$.  Let
$\law_\Fhat$ denote the law of $M_i^{\rm T}X_i$ where $X_i$ is the 
$i$th column of $X$. 

\begin{theorem}
Suppose that each $M_i$ satisfies (\ref{mconds}) and has
$\overline{M_i}=0$. 
Suppose that each $x_i$ satisfies (\ref{Lindeberg}).
Then $\rho_2(\law_\Pbf,\law_{\hat F})\to 0$.
\end{theorem}

The obvious analogue of Theorem \ref{thm4} also holds.

\begin{theorem}
Let $\Delta_0$ be a fixed compact subset of $\Delta$. Suppose $F$ is 
any sequence in $\Delta_0$ and that for each $n$ the $n\times\pi$ 
matrix $X$ has independent rows distributed according
to $F$.  Assume $M$ is an arbitrary sequence of $n\times\pi$ 
matrices 
whose columns $M_i$ each satisfy (\ref{mconds}) and have
$\overline{M_i}=0$.
If $\hat F$ is the empirical distribution function of the rows of $X$
then $\rho_2(\law_F,\law_{\hat F})\to 0$ in probability.
\end{theorem}

It should be noted that the actual null hypothesis plays no role in 
these
theorems.  If the theorems are to be used to deduce that any 
particular
sequence of permutation invariant tests has poor power properties it 
is
necessary that $m_1 = \ldots = m_n$ imply the assertion that the 
null hypothesis is true
and that there be some alternative sequence satisfying the conditions 
of the preceding theorems.

\section{Spacings Statistics}

Suppose $U_1\le\cdots\le U_n$ are the order statistics for a sample 
of
size $n$ from a distribution on the unit interval.  To test the null
hypothesis that this distribution is uniform many authors have
suggested tests based on the sample spacings $D_i=U_i-U_{i-1}$ 
where
we take $U_0=0$ and $U_{n+1}=1$.  Examples of statistics include
Moran's statistic $\sum\log(D_i)$ and Greenwood's statistic $\sum
D_i^2$. See Guttorp and Lockhart (1989) and the references therein 
for
a detailed discussion.  Notice that these statistics are invariant
under permutations of the $D_i$. Also note that the joint distribution
of the $D_i$ is permutation invariant.

Consider a sequence of alternative densities $1+h(x)/n^{1/2}$.
\Cibisov (1961) showed (though his proof seems to rely on 
convergence 
of the permutation moment generating function which does not seem 
to me to follow from the form of the permutation central limit 
theorem
which he cites) under differentiability conditions on $h$ that the
power of any spacings statistic invariant under permutations of the
$D_i$ is asymptotically equal to its level using essentially the
method of proof used above.  We can relax the conditions on $h$
somewhat to achieve the following.

\begin{theorem}
Let $\Delta_0$ be a compact subset of $L_2$, the 
Hilbert space of square integrable functions on the unit interval.
Let $\cal T$ be the family of permutation invariant test functions 
$T(D)$.  As $n\to\infty$,
$$\sup\{|E_h(T)-E_0(T)|: T\in{\cal T}, h\in\Delta_0\}\to 0$$.
\end{theorem}

In Guttorp and Lockhart(1988) it is established that, under the 
conditions 
of the theorem, $\ell(D)$, the log-likelihood ratio, is equal to $\sum
h_i(D_i-1/(n+1)) - \int h^2(x)\, dx/2 +O_P(1)$ for any sequence of
alternatives $h$ converging in $L_2$ where the $h_i$ are suitable 
constants derived from $h$. Using the remark following Lemma 1, 
the fact that 
the joint distribution of the $D_i$ is permutation invariant under the 
null
hypothesis, and the characterization of the spacings as $n+1$ 
independent
exponentials divided by their total, the theorem may be proved by 
following
the argument leading to Theorem 2.

\section{Discussion}

\subsection{Relevance of Contiguity}
All the theorems establish that permutation invariant tests are much 
less
powerful than the Neyman-Pearson likelihood ratio test for 
alternatives
which are sufficiently different from the null that the Neyman-Pearson
test has non-trivial power.  Thus if, in practice, it is suspected which
parameters $m_i$ are the ones most likely to be different from all 
the
others  there will be scope for much more sensitive tests than the
invariant tests.  

On the other hand, some readers will argue that the analysis of 
variance
is often used in situations where no such prior information is 
available.
Such readers, I suspect, will be inclined to argue that this sort of
contiguity calculation is irrelevant to practical people. Some readers 
may feel that a user would need unreasonable amounts of prior 
knowledge to derive a better test than the $F$-test.  Consider the 
situation of the normal example in the first section.  Suppose that the 
$m_i$ can be sorted so that adjacent entries are rather similar.  
Define $h(s) = m_{[ns]}$.  If the sequence of functions are reasonably 
close to some square integrable limit $\eta$ then, without knowing 
$\eta$, we 
can construct a test whose power stays larger than its level if the 
Neyman-Pearson test has the same property.  Specifically consider 
the exponential family example. Let $\gamma_i; i=1,2,\ldots$ be an 
orthogonal basis of $L_2[0,1]$  with each $\gamma_i$ continuous and 
let $\lambda_i$ be any sequence of summable positive constants.  
Define a test statistic of the form 
$T(X) = \sum_i\lambda_i(\sum \gamma_i(j/n) 
X_j)/(\sum_j\gamma_i^2(j/n))^{1/2}$.  If the sequence $h$ converges 
to some $\eta\ne 0$ in $L_2[0,1]$ then the asymptotic power of $T$ 
will be larger than its level.  The test is the analogue of the usual sort 
of quadratic goodness-of-fit test of the Cramer-von~Mises type.

It is worth noting that the calculations compute the power function 
by
averaging over alternative vectors which are a permutation of a 
basic
vector $m$. Another approach to problems with large numbers of 
different
populations (labelled here by the index $i$) is to model the $m_i$
themselves as an iid sequence chosen from some measure $G$.  In 
this case
the null hypothesis is that $G$ is point mass at some unknown value
$\mbar$.  I note that alternative measures $G$ which make the 
resulting
model contiguous to the null make ${\rm Var}_G(m_i)=O(n^{-1/2})$ 
which
means that a typical $m_i$ deviates from $\mbar$ by the  $n^{-1/4}$
discrepancy which arises in our first example and in analysis of 
spacings
tests. In other words when any ordering of the $m_i$ is equally 
likely
vectors $m$ differing from $\mbar$ by the amount we have used 
here are
indistinguishable from the null hypothesis according to this empirical
Bayes model. It is important to note, however, that for this empirical
Bayes model the hypothesis of permutation invariance of the statistic
$T$ is unimportant: if ${\rm Var}_G(m_i)=o(n^{-1/2})$ then every 
test
statistic, permutation invariant or not, has power approaching its 
level.

The proofs hinge rather critically on the exact invariance properties 
of
the statistics considered.  In the Neyman-Scott problem for instance 
if a
single sample size were to differ from all the others the whole 
argument
would come apart.  As long as the sample sizes are bounded the ARE 
of ANOVA is 0 nevertheless, as may be seen by direct calculation 
with the 
alternative non-central $F$-distribution. In the spacings problem of 
section 5  
the sample 2-spacings
defined by $D_i=U_{i+2}-U_i$ still provide tests with non-trivial 
power only at alternatives at the $n^{-1/4}$ distance; the joint 
distribution of these
2-spacings is not permutation invariant and our ideas do not help. 
Our ideas do apply, however, to the non-overlapping statistics of 
Del~Pino(1971).

The definition of ARE offered here may well be challenged since the 
comparison is relative to the Neyman-Pearson test which would not be
used for a composite null versus composite alternative situation. 
Nevertheless there seems to us to be a sharp distinction between procedures
for which our definition yields an ARE of 0 and the quadratic tests
mentioned above whose ARE is then positive.

\subsection{Open Problems and Conjectures}

The results presented here lead to some open problems and obvious 
areas
for further work.  \Hajek's proof of the permutation central limit 
theorem
guarantees convergence of the characteristic function and moments 
up to
order 2 of the variables $\mtPx$.  Our heuristic calculations suggest 
that
a good deal more information could be extracted if the characteristic
function could be replaced by the moment generating function and if
convergence of the moment generating function could be established 
not only for fixed arguments but for arguments growing at rates 
slower 
than $n^{1/2}$.  Such an extension would eliminate, as in the normal 
example, the need for considering only contiguous alternatives.  
Large sample
theory for spacings statistics suggests that the results presented here
hold out to alternatives at such distances.  If, in addition,
approximations were available to the moment generating function 
for
arguments actually growing at the rate $n^{1/2}$ the technique 
might
extend to providing power calculations for alternatives so distant 
that permutation invariant tests have non-trivial limiting power.  
Another possible extension would use Edgeworth type expansions in 
the 
permutation central limit theorem to get approximations for the 
difference 
between the power and the level in the situation, covered by our 
theorems, where
this difference tends to 0.

Consider the exponential family model of section 3 for the special 
case of the normal 
distribution. The problem of testing $m_1 = \cdots = m_n$ is 
invariant under the 
permutation group and under the sub-group of the orthogonal group 
composed of
all orthogonal transformations fixing the vector $\one$ all of whose 
entries are 1. It is 
instructive to compare our results for the two different groups. The 
example illustrates 
the trade-off.  For statistics invariant under larger groups it may be 
easier to
prove the required convergence of the average likelihood ratio; the 
easier proof
is balanced against applying the conclusion to a smaller family of 
statistics.

For statistics invariant under the larger group of orthogonal 
transformations fixing $\one$ we
can modify the argument of section 2 and extend the conclusion 
described in the heuristic
problem of section three to relatively large values for $\norm{m-
\mbar}$. Rather than describe
the details we follow a suggestion made to us by Peter Hooper.  
Suppose $Y=Xb+se$ where $e$
has a multivariate normal distribution with mean vector $m$ and 
variance covariance matrix the
identity and where $s$ is an unknown constant, $X$ is an $n\times 
p$ matrix of regression
covariate values of rank $p$ and $b$ is an unknown $p$ dimensional 
vector. Let $H= X (X^{\rm T}X)^{-1}
X^{\rm T}$ be the hat matrix.
Suppose we wish to test the hypothesis $(I-H)m=0$. 
(If $X=\one$ this is equivalent to the problem mentioned above of 
testing
$m_1 = \cdots = m_n$. 
The problem is identifiable only if $Hm=0$ or equivalently if $X^{\rm 
T}m=0$.)
The problem is invariant under the group $O_X$ of orthogonal 
matrices $\Pbf$
for which $\Pbf X = X$.  

Suppose $T(Y)$ is a family of statistics such that $T(\Pbf Y)=T(Y)$ for 
any $\Pbf$ in $O_X$.
Consider the likelihood ratio of $m$ to $Hm$ (the latter is a point in 
the null).
Following equations (1-3) we are lead to study 
\[\Lbar (Y)=\int L(\Pbf Y) F(d\Pbf ).\]  where now $F$ is Haar 
measure on $O_X$ and 
$L(Y)= \exp(m^{\rm T}(I-H)(Y-Xb) - \norm{(I-H)m}^2/2)$. Since 
$\Pbf X = x$ we see that 
$L(\Pbf Y) = \exp(m^{\rm T}(I-H)\Pbf (Y-Xb) - \norm{(I-H)m}^2/2)$.  
If $\Pbf$ is distributed according
to $F$ and $Z$ is standard multivariate normal then $\Pbf^{\rm T}(I-
H)m/\norm{(I-H)m}$ and $(I-H)
Z/\norm{(I-H)Z}$ have the same distribution. This fact and 
expansions of Bessel functions
show that  $\Lbar \to 1$ in probability provided $\norm{(I-
H)m}=o((n-p)^{1/4})$.

The family of statistics invariant under the group of permutations of 
the entries of $Y$ will
be different than the family invariant under $O_X$.  When $X$ is 
simply a vector which is 
a non-zero multiple of $\one$ the family of statistics invariant under 
the permutation group
is much larger than the family invariant under $O_\one$.  For this 
case we are led to study
the variable $\Lbar = \sum_\Pbf \exp( (m-\mbar)^{\rm T} \Pbf e - 
\norm{m-\mbar}^2/2)/n!$.  We 
find that $\E_\mbar(\Lbar)=1$ and $\Var_\mbar(\Lbar)= 
\sum_\Pbf \exp( (m-\mbar)^{\rm T} \Pbf
(m-\mbar))/n! - 1$.  Just when this variance goes to 0 depends on 
extending the permutation
central limit theorem to give convergence of moment generating 
functions.  Since the random
variable $(m-\mbar)^{\rm T} \Pbf (m-\mbar)$ has mean 0 and 
variance $\norm{m-\mbar}^4/n$ we
are again led to the heuristic rate $\norm{m-\mbar}=o(n^{1/4})$.  
However, by taking $m$ to
have exactly one non-zero entry it is not too hard to check that this 
heuristic calculation
cannot be made rigorous without further conditions on $m$ to 
control the largest entries.

Finally, if $X$ is not a scale multiple of $\one$ the problem is not 
invariant under the
permutation group.  Is there some natural extension of our 
techniques to this context for a group smaller than $O_X$?

\appendix
\section*{Appendix}

\subsection*{Proof of Theorem \ref{thm4}}
 
We prove below (cf Shorack and Wellner, p 63 their formula 5, 
except that there the distribution $F$ does not depend on $n$) that
\begin{equation}
\label{empirical}
\rho_2(F,\Fhat)\to 0
\end{equation}
in probability. If Theorem~\ref{thm4} were false then from any counterexample
sequence we could extract a subsequence which is a counterexample 
and along which the convergence in (\ref{empirical}) is almost sure.
The theorem then follows from the assertion that 
\begin{equation}
\label{unifint}
\rho_2(F,G)\to 0 \qquad {\rm implies }\qquad \rho_2(\law_F,\law_G)\to 0
\end{equation} 
whenever $F$ is any sequence of distributions with compact
closure in $\Delta$. Assertion (\ref{unifint}) is a consequence
of Lemma 1 of Guttorp and Lockhart (1988).

To prove (\ref{empirical}) we may assume without loss, in view of the
compactness of $\Delta_0$ that $F\dconv\Phi$ for some $\Phi$. 
Elementary moment calculations assure that $\Fhat(\tau)$ converges in 
probability to
$\Phi(\tau)$ for each $\tau$ which is a continuity point of $\Phi$. 
This guarantees that $\rho_0(\Fhat,\Phi)\to 0$ in probability.
We need only
show that $\E_\Fhat(X^2)\to\E_\Phi(X^2)$.  But $\E_\Fhat(X^2)=\sum
X_i^2/n$.  The triangular array version of the law of large
numbers given in Lemma 2 of Guttorp and Lockhart(1988) shows that 
$\sum X_i^2/n - \E_F(X^2)\to 0$ in probability.  Since $F\dconv\Phi$ 
implies that $\E_F(X^2)\to \E_\Phi(X^2)$ we are done.

\subsection*{Proof of Theorem \ref{nonexpthm}}

It remains to choose a sequence $k=k(n)$ in such a way that we can
check (\ref{step3}).  In view of permutation invariance we may 
assume without loss that $|m_1-\mbar| \ge \cdots \ge |m_n-\mbar|$.
Define matrices $C$ and $D$ by setting $C(i,j) = \phi(X_j,m_i)-
\phi(X_j,\mbar)
- \E_\mbar ( \phi(X_j,m_i)-\phi(X_j,\mbar) )$  and $D(i,j) 
= (m_i - \mbar) \phi_1(X_j;\mbar)$. We will eventually choose a 
sequence $k$ and
put $B(i,j)=C(i,j)$ for $i\le k$ and $B(i,j)=D(i,j)$ for $i > k $. 
Note that $\sum_i B(i,i) $ is simply $S^*(X)$.  

If $\Pbf$ is a random 
permutation matrix then in row $i$ there is precisely 1 non-zero entry;
let $J_i$ be the column where this entry occurs.  Then $S^*(\Pbf X) =
\sum_i B(i,J_i)$.  The variables $J_1,\ldots,J_n$ are a random permutation
of the set $\{1,\ldots,n\}$.  As in the proof of Theorem 2 let $J_1^*,\ldots,
J_n^*$ be a set of independent random variables uniformly distributed on
$\{1,\ldots,n\}$. We will show that for each fixed $\kappa$ 
\begin{equation}
\label{slow1}
\rho_2\left(\law(\sum_{i\le\kappa}C(i,J_i^*)|X),
\law(\sum_{i\le\kappa}C(i,i))\right)\to 0
\end{equation}
and
\begin{equation}
\label{slow2}
\rho_2\left(\law(\sum_{i\le\kappa}C(i,J_i^*)|X),
\law(\sum_{i\le\kappa}C(i,J_i)|X)\right)\to 0
\end{equation}
in probability.  We will also show that for any sequence $k$ tending 
to $\infty$ with $k^2=o(n)$ we have
\begin{equation}
\label{slow3}
\rho_2\left(\law(\sum_{i > k}(m_i-\mbar)\phi_1(X(J_i^*),\mbar)|X), 
\law(\sum_{i > k}D(i,i))\right)\to 0 
\end{equation} 
in probability.  There is then a single sequence $k$ tending to 
infinity so slowly that (\ref{slow1}) and (\ref{slow2}) hold with 
$\kappa$ replaced by $k$ and so that (\ref{slow3}) holds.  We use 
this sequence $k$ to define $B$.

We will then show that 
\begin{equation}
\label{boot}
\rho_2\left(\law (\sum_i B(i,J_i^*) | X ), \law(\sum_i B(i,i))\right) \to 0
\end{equation}
in probability and that 
\begin{equation}
\label{srs}
\rho_2\left(\law(\sum_i B(i,J_i)| X),\law (\sum_i B(i,J_i^*)|X)\right) \to 0
\end{equation}
in probability.  These
two are enough to imply (\ref{step3}) as in Corollary 2 and the 
obvious (but unstated) corresponding corollary to Theorem \ref{thm3}. 

\vskip12pt
\leftline{\em Proof of (\ref{slow1})}

For each fixed $i$ we may apply (\ref{empirical}) with
the vector $(X_1,\ldots,X_n)$ replaced by $(C(i,1),\ldots,C(i,n))$ to 
conclude that 
$$\rho_2(\law(C(i,J_i^*)|X),\law(C(i,i)))\to 0$$
in probability; the condition imposed on $F$ leading to (\ref{empirical})
is implied by (\ref{scoreassump}).  Use the independence properties 
to conclude
$$\rho_2(\law(C(1,J_1^*),\ldots,C(\kappa,J_\kappa^*)|X),
\law(C(1,1),\ldots,C(\kappa,\kappa)))\to 0$$ for each fixed 
$\kappa$.
Assertion (\ref{slow1}) follows. 

\vskip12pt
\leftline{\em Proof of (\ref{slow2})}

For each fixed $i$ we have $\law(C(i,J_i^*)|X)=\law(C(i,J_i)|X)$.  
Furthermore it is possible to construct $J$ and $J^*$ in such a way
that for each fixed $\kappa$ we have $P(J_i=J_i^*; 1\le 
i\le\kappa)\to 1$.
This establishes (\ref{slow2}).

\vskip12pt
\leftline{\em Proof of (\ref{slow3})}

Let $\mbar_{-k}=\sum_{i > k}m_i/(n-k)$ and let $U$ be the vector with $i$th entry $\phi_1(X_i,\mbar)$.  Arguing as in the proof of 
Theorem \ref{thm4} (see \ref{unifint} above) we see that  
\begin{equation}
\label{law1}
\rho_2\left(\law(\sum_{i > k}(m_i-\mbar_{-k})U(J_i^*)|X),
\law(\sum_{i > k}(m_i-\mbar_{-k})U_i)\right)\to 0
\end{equation}
in probability. We need to replace $\mbar_{-k}$ by $\mbar$ in order
to verify (\ref{slow3}).  Elementary algebra shows 
$\mbar_{-k}-\mbar=O(k/(n-k))$.  Temporarily let 
$$T_1 = \sum_{i >
k}(m_i-\mbar_{-k})U(J_i^*) - \sum_{i >
k}(m_i-\mbar)U(J_i^*),$$  and
$$T_2 = \sum_{i > k}(m_i-\mbar_{-k})U_i- \sum_{i >
k}(m_i-\mbar)U_i $$
In view of (\ref{scoreassump}) we see that $\Var(U_1)=O(1)$.  
Hence
\begin{equation}
\label{var1}
\Var(T_2)  = 
(\mbar_{-k}-\mbar)^2(n-k)\Var(U_1)\to 0.
\end{equation}
Since the central limit theorem shows that the sequence 
$\law(\sum_{i > k}D(i,i))$ has compact closure in $\Delta$ we may 
use (\ref{var1}) to show that 
\begin{equation}
\label{law2}
\rho_2\left(\law(\sum_{i > k}(m_i-\mbar_{-k}) U_i), 
\law(\sum_{i > k}D(i,i))\right) \to 0.
\end{equation}

Next
\begin{equation}
\label{var2}
\Var(T_1 |X) =  (\mbar_{-k}-\mbar)^2(n-k) 
\Var(U(J_1^*)|X).
\end{equation}
Since 
$$\Var(U(J_1^*)|X)= 
\sum U^2_i/n - (\sum U_i/n)^2$$
we may apply
the triangular array law of large numbers given in Guttorp and 
Lockhart(1988, Lemma 2) to conclude that 
$\Var(U(J_1^*)|X)= O_P(1)$. 
Since $k^2=o(n)$  we see that the right hand side of (\ref{var2}) 
tends to 0 in probability.  Hence
\begin{equation}
\label{law3}
\rho_2\left(\law(\sum_{i > k}(m_i-\mbar_{-k})U(J_i^*)|X),
\law(\sum_{i > k}(m_i-\mbar)U(J_i^*)|X)\right)\to 0
\end{equation}
in probability. Assembling (\ref{law1}), (\ref{law2}) and (\ref{law3}) 
we have 
established (\ref{slow3}).

\vskip12pt
\leftline{\em Proof of (\ref{boot})}

Given $X$,  the variables $\sum_{i\le k}C(i,J_i^*)$ and 
$\sum_{i > k}(m_i-\mbar)\phi_1(X(J_i^*),\mbar)$ are independent.  
Similarly $\sum_{i\le k}C(i,i)$ and $\sum_{i > k}D(i,i)$ are 
independent.  Statement (\ref{boot}) then follows from (\ref{slow1}) 
and (\ref{slow3}).

\vskip12pt
\leftline{\em Proof of (\ref{srs})}

To deal with (\ref{srs}) we must cope with the lack of independence 
among the  $J_i$.  A random permutation of
$\{1,\ldots,n\}$ can be generated as follows.  Pick $J_1,\ldots,J_k$ a
simple random sample of size $k$ from $\{1,\ldots,n\}.$ Let $\Pbf_{-k}$,
independent of $J_1,\ldots,J_k$, be a random
permutation of $\{1,\ldots,n-k\}.$  Then $\sum B(i,J_i)$ has the same
distribution (given $X$) as
\[\sum_{i=1}^kC(i,J_i) +
(m_{-k}-\mbar_{-k})^{\rm T}\Pbf_{-k}U_{-k}\] where the subscript 
$-k$ on $m$ denotes deletion of $m_1,\ldots,m_k$ while that on $U$ denotes 
deletion of the entries $U(J_1),\ldots,U(J_k)$.

Let $Z_0$ denote a random variable, independent of $J$ and $J^*$
whose conditional distribution given $X$ is normal with mean 0 and 
variance $\sum_{i>k}(m_i-\mbar)^2S$ where $S$ is the sample variance of 
the $U_j$'s. Statement (\ref{srs}) is a consequence of the following 3
assertions:
\begin{equation}
\label{f1}
\rho_2\left(\law(\sum B(i,J_i^*)|X),\law(\sum_{i\le k}C(i,J_i^*)+Z_0|X)\right)
\to 0
\end{equation}
in probability,
\begin{equation}
\label{f2}
\rho_2\left((\law(\sum_{i\le k}C(i,J_i)+Z_0|X), \law(\sum_{i\le 
k}C(i,J_i^*)+Z_0|X)\right)
\to 0
\end{equation}
in probability,
and
\begin{equation}
\label{f3}
\rho_2\left(\law(\sum_{i\le k}C(i,J_i)+Z_0|X), 
\law(\sum_iB(i,J_i)|X)\right)
\to 0
\end{equation}
in probability.

Condition (\ref{f2}) follows from (\ref{slow2}), the conditional
independence of $Z_0$ and $J,J^*$ and the fact that the conditional
variance of $Z_0$ is bounded.  Condition (\ref{f1}) is implicit in the 
proof of
(\ref{boot}) after noting that the variance $S$ of the entries in $U$
is negligibly different from $\iota(\mbar)$. 

It remains to establish (\ref{f3}).  We will condition on 
$(J_1,\ldots,J_k)$ as well as $X$ and apply the Permutation
Central Limit Theorem.  The application of the conditions of that 
theorem is a bit delicate since the conditions will only be shown 
to hold in probability. We present the argument in the form of a 
technical lemma.  

\begin{lemma}
Suppose $(W_1,W_2,W_3)$ is a sequence of random variables.  
Suppose $g$ and $h$ are sequences of measurable functions defined 
on the range spaces of $(W_1,W_2)$ and $(W_1,W_2,W_3)$.  Let 
$\zeta_1,\zeta_2$ be two independent real valued random 
variables. Suppose that there are functions $f_i$ for $i=0,1,\ldots$ 
(also indexed as usual by the hidden index $n$) such that
\begin{equation}
\label{aplemcond1}
f_0(w_1)\to 0\qquad {\rm implies }\qquad g(w_1,W_2)\dconv 
\zeta_1
\end{equation}
and 
\begin{equation}
\label{aplemcond2}
f_1(w_1,w_2)\to 0\qquad i=1,2,\ldots\qquad{\rm implies }\qquad 
h(w_1,w_2,W_3)\dconv \zeta_2\, .
\end{equation}
If $f_0(W_1)\to 0$ in probability and $f_i(W_1,W_2)\to 0$ in 
probability then 
$$\law(g(W_1,W_2)+h(W_1,W_2,W_3)|W_1) \Rightarrow 
\zeta_1+\zeta_2$$
in probability.  If in addition
\begin{equation}
\label{aplemcond3}
\E(g(W_1,W_2)+h(W_1,W_2,W_3)|W_1) \to \E(\zeta_1+\zeta_2)
\end{equation}
in probability,
\begin{equation}
\label{aplemcond4}
\Var(g(W_1,W_2)|W_1) \to \Var(\zeta_1)
\end{equation}
in probability and 
\begin{equation}
\label{aplemcond5}
\Var(h(W_1,W_2,W_3)|W_1) \to \Var(\zeta_2)
\end{equation}
in probability then 
$$\law(g(W_1,W_2)+h(W_1,W_2,W_3)|W_1) \dconv 
\zeta_1+\zeta_2$$
in probability.
\end{lemma}

The lemma is to be applied with $W_1 = X$, with $W_2 = 
(J_1,\ldots,J_k)$ and with $W_3 = \Pbf_{-k}$.  Conclusion (\ref{f3})
can be reduced to the form given by a compactness argument.  The sequence 
of laws of $\sum_{i\le k} C(i,i)$ has compact closure in $\Delta$ in 
view of (\ref{scoreassump}).  Then apply (\ref{slow1}) and (\ref{slow2}) 
to conclude that the sequence of laws of $\sum_{i\le k}C(i,J_i)$ also 
has compact closure.  Let the distribution of $\zeta_1$ be any limit 
point of this sequence of laws.  The 
random variable $\zeta_2 $ will be the normal limit in distribution of
$(m_{-k}-\mbar)^{\rm T}\Pbf_{-k}U_{-k}$. 

In order to apply the lemma we must give the conditions of 
the Permutation Central Limit
Theorem in a form in which we have only a countable family of 
convergences, as required in (\ref{aplemcond2}), to check.
Note that (\ref{Lindeberg})  in Theorem \ref{thm3} 
can be replaced by the assertion that
there is a sequence $\tau_1,\ldots$ of real
numbers increasing to $\infty$ such that 
\begin{equation}
\label{pclt2}
 \max(n^{-1}\sum(x_i-\xbar)^21(|x_i-\xbar| >
\tau_j),\Psi(\tau_j))-\Psi(\tau_j)\to 0
\end{equation}
for each $j\ge 1$.

We now apply Theorem \ref{thm3} with $x$ replaced by $U_{-k}$, with $m$ 
replaced by $m_{-k}-\mbar$ and with $\Pbf$ replaced by $\Pbf_{-k}$. 
It is easy to check by conditioning on $J$ that $\E(U_{-k}) = 0$ and 
$\Var(\overline{U_{-k}})= \iota(\mbar)/(n-k)$.  Hence$(n-k) (\overline{m_{-k}}-\mbar)
\overline{U_{-k}}\to 0$ in probability.  Set
$$ T_3 = \sum(U_{-k,i}-\overline{U_{-k}})^2 
1(| U_{-k,i}-\overline{U_{-k}} | > t)/(n-k)$$
Then
\begin{eqnarray*}
T_3 & \le& 2 \sum U^2_{-k,i} 1(| U_{-k,i}-\overline{U_{-k}} | > t)/(n-k)\cr &
& +2 \overline{U_{-k}}^2 \sum 1(| U_{-k,i}-\overline{U_{-k}} | > t)/(n-k).
\end{eqnarray*}
The second term on the right is $O_P(1/(n-k))=o_P(1)$.
Since 
$$1(| U_{-k,i}-\overline{U_{-k}} | > t) \le 
1( |U_{-k,i} | > t/2) + 1( | \overline{U_{-k}}| > t/2)$$ 
we see that 
$$
T_3  \le  2 (n-k)^{-1} \sum_i U^2_i1(|U_i| >t/2) +o_P(1)\, .
$$
Take $\Psi(t)= 2\sup\{\E_m(\phi_1(X_1,m)1(|\phi_1(X_1,m)| > t/2) ;
m\in\Theta_0\}$ and apply Lemma 2 of Guttorp and Lockhart together with 
(\ref{scoreassump}) to check that (\ref{pclt2}) holds.  To finish the
proof of (\ref{f3}) we need to check convergence of second moments as in 
(\ref{aplemcond3}), (\ref{aplemcond4}) and (\ref{aplemcond5}).  This can
be done using (\ref{slow1}), (\ref{slow2}) and direct calculation of the
conditional mean and variance given $X$ of $\sum_{i>k}D(i,J_i)$. 
Theorem (\ref{nonexpthm}) follows.

The technical lemma itself may be proved as follows.  From any
counterexample sequence we may extract by a diagonalization 
argument
a subsequence which is still a counterexample and for which 
$f_0(W_1)\to 0
$ almost surely and $f_j(W_1,W_2)\to 0$ almost surely 
for each $j \ge 1$.  For
any sample sequence for which all these convergences occur we have
$\law(h(W_1,W_2,W_3)|W_1,W_2)\Rightarrow \zeta_2$ and 
$\law(g(W_1,W_2)| W_1)\Rightarrow \zeta_1$. Evaluation of the
 conditional characteristic function of 
$g(W_1,W_2)+h(W_1,W_2,W_3)$ given $W_1$ by further 
conditioning on $W_2$ yields the convergence in distribution 
asserted in the lemma.  The remaining conclusions concerning 
moments are more elementary analogues of the same idea.

\section*{References}

\myref{Abramowitz, M. and Stegun, I. A.
(1965).
{\it Handbook of Mathematical Functions}. New York: Dover.}

\myref{Billingsley, Patrick
(1968).
{\it Convergence of Probability Measures}. New York: Wiley.}

\myref{\Cibisov, D.M.
(1961).
On the tests of fit based on sample spacings. {\it Teor.
Verojatnos. i Primenen.} {\bf 6}, 354--8.}

\myref{Del Pino, G.E.
(1979).
On the asymptotic distribution of k-spacings with applications to
goodness-of-fit tests. {\it Ann. Statist.}, {\bf 7}, 1058-1065.}

\myref{Guttorp, P. and Lockhart, R. A.
(1988).
On the asymptotic distribution of quadratic forms in uniform order
statistics. {\it Ann. Statist.}, {\bf 16}, 433-449.}

\myref{H{\'a}jek, J. 
(1961). 
Some estensions of the Wald-Wolfowitz-Noether Theorem.
{\it Ann. Math. Statist.}, {\bf 32}, 506-523.}

\myref{H{\'a}jek, J. and {\v S}id{\'a}k, Z. (1967). {\it Theory of Rank
Tests}.  Academic Press: New York.}

\myref{Shorack, G. R. and Wellner, J. A. (1986).  {\it Empirical
Processes with Applications to Statistics}. New York: Wiley.}

\end{document}